\documentclass[12pt]{article}

\usepackage{amsmath}
\usepackage{amsfonts}

\newtheorem{theo}{Theorem}
\newtheorem{coro}{Corollary}

\newtheorem{lemm}{Lemma}
\newtheorem{defn}{Definition}

\def\E{\mathbb{E}}
\def\PP{\mathbb{P}}
\def\0{{\bf 0}}
\def\Z{\mathbb{Z}}
\def\R{\mathbb{R}}

\def\rc{{\rm c}}

\renewcommand{\E}{\mathbb E \,}

\newcommand{\tod}{\stackrel{{\cal D}}{\longrightarrow}}

\newcommand{\eqco}{\setcounter{equation}{0}}
\newcommand{\thco}{\setcounter{theo}{0}}
\newcommand{\prco}{\setcounter{prop}{0}}
\newcommand{\laco}{\setcounter{lemm}{0}}
\newcommand{\coco}{\setcounter{coro}{0}}
\newcommand{\cjco}{\setcounter{conj}{0}}

\newcommand{\deco}{\setcounter{defn}{0}}
\newcommand{\allco}{\eqco  \thco \prco \laco \coco \cjco \deco}
\newcommand{\qed}{\hfill{\rule[-.2mm]{3mm}{3mm}}}

\newcommand{\Po}{{\cal P}}
\newcommand{\X}{{\cal X}}
\def\b{{\beta}}

\renewcommand{\H}{{\cal H}}
\renewcommand{\P}{{\cal P}}

\newcommand{\Cov}{{\rm Cov}}

\newcommand{\Var}{{\rm Var}}

\newcommand{\XX}{{\cal X}}
\newcommand{\MM}{{\cal M}}

\newcommand{\F}{{\cal F}}

\newcommand{\NN}{{\cal N}}

\newcommand{\SIG}{{\rm SIG}}

\def\bdm{\begin{displaymath}}
\newcommand{\edm}{\end{displaymath}}
\def\benu{\begin{enumerate}}
\def\eenu{\end{enumerate}}
\def\beqn{\begin{equation}}
\def\eeqn{\end{equation}}
\def\be{\begin{equation}}
\def\ee{\end{equation}}
\def\bea{\begin{eqnarray}}
\def\eea{\end{eqnarray}}
\newcommand{\bean}{\begin{eqnarray*}}
\newcommand{\eean}{\end{eqnarray*}}
\newcommand{\bear}{\begin{eqnarray}}
\newcommand{\eear}{\end{eqnarray}}

\newcommand{\1}{{\bf 1}}

\renewcommand{\epsilon}{\varepsilon}

\newcommand{\stau}{{\nu}}

\def\R{\mathbb{R}}

\def\al{{\alpha}}
\def\de{{\delta}}

\def\qed{\hfill\hbox{${\vcenter{\vbox{
    \hrule height 0.4pt\hbox{\vrule width 0.4pt height 6pt
    \kern5pt\vrule width 0.4pt}\hrule height 0.4pt}}}$}}

\def\la{{\lambda}}
\def\ka{{\kappa}}

\begin{document}

\title{\bf Normal Approximation in Geometric Probability}

\author{Mathew D. Penrose$^1$ and J. E. Yukich$^{2}$ \\
\\
{\normalsize{\em University of Bath and Lehigh University}} }
%\date{Revised version, August 2004}
\maketitle

\footnotetext{$~^1$ Department of Mathematical Sciences, University
of Bath, Bath BA1 7AY, United Kingdom: {\texttt m.d.penrose@bath.ac.uk} }

\footnotetext{$~^2$ Department of Mathematics, Lehigh University,
Bethlehem PA 18015: {\texttt joseph.yukich@lehigh.edu} }

%\footnotetext{$~^1$ Research supported in part by the Isaac
%Newton Institute for Mathematical Sciences, Cambridge}

\footnotetext{$~^2$ Research supported in part by NSA grant
MDA904-01-1-0029 and NSF grant DMS-0203720}

\begin{abstract}
We use Stein's method to obtain bounds on the rate of convergence
for a class of statistics in geometric probability obtained as a
sum of contributions from Poisson points which are exponentially
stabilizing, i.e. locally determined in a certain sense. Examples
include statistics such as total edge length and total number of
edges of graphs in computational geometry and the total number of
particles accepted in random sequential packing models. These
rates also apply to the 1-dimensional marginals of the random
measures associated with these statistics.
\end{abstract}

\newpage

\section{Introduction}
\allco

In the study of limit theorems for functionals on Poisson or
binomial spatial point processes, the notion of {\em
stabilization} has recently proved to be a useful unifying concept
\cite{BY2,PY1,PY4}. Laws of large numbers and central limit
theorems can be proved in the general setting of functionals
satisfying an abstract `stabilization' property whereby the
insertion of a point into a Poisson process has only a local
effect in some sense. These results can then be applied to deduce
limit laws  for a great variety of particular functionals,
including those concerned with the minimal spanning tree, the
nearest neighbor graph, Voronoi and Delaunay graphs, packing, and
germ-grain models.

Several different techniques are available for
proving general central limit theorems for stabilizing
functionals. These include a martingale approach
\cite{PY1} and  a method of moments
 \cite{BY2}.
In the present work, we revisit a third technique
for proving central limit theorems for stabilizing functionals
on Poisson point processes, %in geometric probability,
which was introduced by Avram and Bertsimas \cite{AB}.
This method is based on the  normal approximation
of sums of random variables which are `mostly independent
of one another' in a sense made precise via
 dependency graphs, which in turn is proved  via
Stein's method \cite{St}.  It has the advantage of providing
 explicit error bounds and rates of convergence.

We extend the work of Avram and Bertsimas in several directions.
First, whereas in \cite{AB} attention was restricted to certain
particular functionals, here we derive a general result holding
for arbitrary functionals satisfying a stabilization condition
which can then be checked rather easily for many special cases.
Second, we consider
non-uniform point process intensities
and do not require the functionals to be translation invariant.
Third, we improve on the
rates of convergence in \cite{AB} by making use of the recent
refinement by Chen and Shao \cite{CS} of previous normal
approximation results for sums of `mostly independent' variables.
Finally, we apply the methods not only to random variables
obtained by summing some quantity over Poisson points,
but to the associated {\em random point measures}, thereby
recovering many of the results of Baryshnikov and Yukich
\cite{BY2} on convergence of these measures,
and without requiring higher order moment calculations.
We add to  \cite{BY2} by
 providing  information about the rate of convergence,
and relaxing the continuity conditions
required in \cite{BY2} for
test functions and point process intensities.

A brief comparison between the methods of deriving central limit
theorems for functionals of spatial point processes is warranted.
  Only the dependency graph method used
here, to date, has yielded error bounds and rates of convergence.
On the other hand, our method requires bounds on the tail of the
`radius of stabilization' (i.e., on the range of the local effect
of an inserted point). The martingale method, in contrast,
requires only that this radius be almost surely finite, and for
this reason is applicable to some examples such as those concerned
with the minimal spanning tree,  for which no tail bounds are
known and which therefore lie beyond the scope of the present
work. The moment method \cite{BY2}  and martingale method \cite{Pe2},
 unlike the dependency graph method,
provide information about the variance of the Gaussian limits. The
moment method  has also been used \cite{BESY} to establish
%describes the gaussian and
 moderate scale limit behavior of
functionals of spatial point processes.
%in terms of the underlying
%density of points \cite{BY2, BESY}. (does martingale method also
%yield limits like those in \cite{BY2}?).
 Whereas the moment method
requires exponential tail bounds for the radius of stabilization,
%The method of moments \cite{BY2} does not provide error
%bounds, but otherwise seems to have a similar scope to the methods
%used here, in that it also requires tail bounds on the radius of
%stabilization.
one of our central limit theorems (Theorem \ref{polyth})
 requires only that
this tail $\tau(t)$ decay as a (large) negative power of $t$.
%whereas the central limit theorem in \cite{BY2} requires $\tau(t)$
%to decay exponentially in $t$.

With regard to ease of use in applications, the dependency
graph method and method of moments require checking
tail bounds for the radius of stabilization, which is
usually straightforward where possible at all.
The method of moments requires a more complicated
(though checkable) version of the bounded moments
condition (\ref{mom}) below (see \cite{BY2}). The
dependency graph method requires some separate calculation
of variances if one wishes  to identify explicitly the
variance of the limiting normal variable. The martingale
method requires the checking of slightly more subtle
versions of the stabilization conditions needed here \cite{Pe2,PY1}.

\section{General results}
\label{secgr}

Let $d \geq 1$ be an integer.
For the sake of generality, we consider {\em marked} point
processes in $\R^d$. Let $(\MM,\F_\MM,\PP_\MM)$ be a probability
space (the {\em mark space}).
 Let $\xi((x,s);\X)$ be a measurable
$\R$-valued function  defined for all pairs $((x,s),\X)$, where
$\X \subset \R^d \times \MM$
 is finite and where $(x,s) \in \X$ (so $x \in \R^d$ and $s \in \MM$).
When $(x,s) \in (\R^d \times \MM) \setminus \X$, we
 abbreviate notation and write $\xi((x,s); \X)$
instead of $\xi((x,s); \X \cup \{(x,s)\})$.

Given $\X \subset \R^d \times \MM$,  $a > 0$ and $y \in \R^d$, we let $a\X:=
\{(ax,t): (x,t) \in \X\}$ and $y+ \X := \{(y+x,t): (x,t) \in \X \}$;
in other words, translation and scalar multiplication on $\R^d \times
\MM$ act only on the first component.
 For all $\la > 0$ let
$$
\xi_{\la}((x,s); \X) := \xi ((x,s); x+ \la^{1/d}(-x +\X) ).
$$

We say $\xi$ is {\em translation invariant}  if $\xi((y+ x ,s); y+
\X ) = \xi((x,s); \X)  $ for all $y \in \R^d$, all $(x,s)\in \R^d
\times \MM$ and all finite $\XX \subset \R^d \times \MM$. When
$\xi$ is
 translation invariant,  the functional $\xi_\la$  simplifies to
$\xi_{\la}((x,s); \X) = \xi ((\la^{1/d}x,s); \la^{1/d}\X )$.

Let $\ka$ be a probability  density function
 on $\R^d$ with compact support $A \subset \R^d$.
For all
$\la > 0$, let
%$\P_{\la}:= \P_{\la \ka}$
$\P_{\la}$
 denote a Poisson point process
 in $\R^d \times \MM$
 with intensity measure $(\la \ka (x) dx) \times \PP_\MM(ds)$.
 We shall assume throughout that $\ka$ is bounded with supremum
denoted $\|\kappa\|_\infty$.

Let $(A_\la, \la \geq 1)$ be a family of Borel subsets  of $A$.
The simplest case, with
 $A_\la = A$ for all $\la$, covers all examples considered
here; we envisage possibly using the general case in future work.

The following notion of exponential stabilization,  adapted from
 \cite{BY2}, plays a central role in all that follows.
 For $x \in \R^d$ and $r > 0$, let $B_r(x)$ denote the Euclidean ball
centered at $x$ of radius $r$. Let $U$ denote a random element of
$\MM$ with distribution $\PP_\MM$, independent of $\Po_\la$.

\begin{defn}\label{stab}  $\xi$ is {\em
exponentially stabilizing} with respect to $\ka$ and $(A_\la)_{\la
\geq 1}$
 if for all $\la \geq 1 $ and all
 $x \in  A_\la$,
there exists an a.s. finite random variable $R:=R(x,\la)$  ({\em a
radius of stabilization} for $\xi$ at $x$) such that for all
finite ${\X } \subset (A \  \setminus  B_{\la^{-1/d} R} (x))
\times \MM$, we have \be \label{stab1}
 \xi_\la \left((x,U);  [ {\cal P}_{\la}  \cap
 (B_{\la^{-1/d} R}(x) \times \MM) ]
  \cup \X \right)
=
 \xi_\la \left((x,U);   {\cal P}_{\la}  \cap
 (B_{\la^{-1/d} R}(x) \times \MM)
   \right),
\ee
and  moreover
the tail probability $\tau(t)$ defined for $t >0$ by
\bea
\tau(t) := \sup_{\la \geq 1,\ x \in A_\la}
P[R(x,\la) > t]
\label{taudef}
\eea
satisfies
\bea
\limsup_{ t \to \infty} t^{-1} \log \tau(t) < 0.
%   \leq C \exp (- t/C).
\label{tau1}
\eea
For $\gamma >0$,
we say $\xi $ is {\em polynomially stabilizing of order} $\gamma$
if  the above conditions hold with (\ref{tau1}) replaced by
the condition $\limsup_{t \to \infty} t^\gamma \tau(t) <\infty$.
\end{defn}

Condition (\ref{stab1}) may be cast in a more transparent form as
follows. Each point of $\X$ is a pair $(x,U)$, with
$x \in \R^d$ and $U \in \MM$, but for notational convenience we can view it
as a point $x$ in $\R^d$ carrying a mark $U:=U_x$.
Then we can view $\X$ as a point set in $\R^d$ with each
point carrying a mark in $\MM$. With this interpretation,
(\ref{stab1}) stipulates that for all finite (marked) ${\X }
\subset A \ \setminus  B_{\la^{-1/d} R} (x) $, we have
 \be
  \label{stab2}
 \xi_\la \left(x;   ({\cal P}_{\la}  \cap B_{\la^{-1/d} R}(x) )
  \cup \X \right)
=
 \xi_\la \left(x;  {\cal P}_{\la}  \cap
 B_{\la^{-1/d} R}(x)  \right).
\ee
 Roughly speaking, $R:=R(x,\la)$ is a radius of stabilization
if the value of $\xi_\la( x; {\cal P}_{\la})$ is
unaffected by changes to the points outside $B_{\la^{-1/d}R}(x)$.

Functionals of spatial point processes often satisfy exponential
stabilization (\ref{stab1}) (or (\ref{stab2})); here is an example.
 Suppose $\MM = [0,1]$ and $\PP_M$ is the uniform distribution on $[0,1]$.
Suppose that $A$ is convex or polyhedral, and $\ka$ is bounded away
 from zero on $A$. Suppose a measurable function $(q(x),x\in A)$
is specified, taking values in $[0,1]$.  Adopting the
 conventions of the preceding paragraph, for a marked point set
 $\X \subset \R^d$ let us denote each point $x \in \X$ as `red' if
$U_x \leq q(x)$ and as `green' if $U_x > q(x)$. Let $\xi(x; \X)$
take the value 0  if  the nearest neighbor of $x$ in $\X$ has the
same color as $x$, and take the value 1 otherwise. Note that $\xi$
is {\em not} translation invariant in this example, unless
$q(\cdot)$ is constant. For $x \in A$ let $R:= R(x, \la)$ denote
the distance between $\la^{1/d}x$ and its nearest neighbor in
$\la^{1/d} \P_\la$. Then stabilization (\ref{stab2}) holds because
points  lying outside $B_{\la^{-1/d}R}( x)$ will not change the
value of $\xi_\la(x; \P_\la)$, and it is easy to see that $R$ has
exponentially decaying tails. This  example is relevant to the
multivariate two-sample test described by Henze \cite{He}. See
Section \ref{secapp} for  further   examples.

\begin{defn} \label{momdef} $\xi$ has a moment of order $p > 0$
(with respect to $\kappa$ and $(A_\la)_{\la \geq 1}$) if
\begin{equation} \label{mom}
  \sup_{\lambda \geq 1, \ x \in
A_\la} \E [ | \xi_{\la}((x,U); {\cal P}_{\la } ) |^p ]
< \infty.
\end{equation}
\end{defn}

For $\la >0$, define the random weighted point measure
$\mu_\la^\xi$ on $\R^d$ by
$$
{\mu}_{\la }^{\xi}:= \sum_{(x,s) \in {\cal P}_{\la} \cap (A_\la
\times \MM)} \xi_{\la}((x,s);{\cal P}_{\la }) \delta_x
$$
and the centered  version $\overline{\mu}^{\xi}_{\la } :=
{\mu}_{\la }^{\xi} - \E[{\mu}_{\la }^{\xi}]$.

Let $B(A)$ denote
the set of bounded Borel-measurable functions on $A$.
Given $f \in B(A)$, let
$\langle f, {\mu}^{\xi}_{\la } \rangle := \int_A f d
{\mu}^{\xi}_{\la }$
and
$\langle f, \overline{\mu}^{\xi}_{\la } \rangle := \int_A f d
\overline{\mu}^{\xi}_{\la }$.

Let $\Phi$ denote the distribution function of the standard
normal.  Our main result
is a normal approximation result for
$\langle f, \overline{\mu}^{\xi}_{\la } \rangle $, suitably scaled.

\begin{theo}\label{mainthm1}
Suppose $\|\kappa\|_\infty < \infty$. Suppose that
  $\xi$ is exponentially stabilizing and
satisfies the moments condition (\ref{mom}) for some $p >2$.
Let $q \in (2,3]$ with $q < p$.
Let $f \in B(A)$ and put $T_\la :=
\langle f, \mu^{\xi}_{\la } \rangle$.
  Then there exists
 a finite constant $C$ depending on $d, \xi$,  $\ka$, $p$, $q$ and $f$,
such that for all $\la \geq 2$,
 \be
\label{rates2} \sup_{t \in \R}  \left| P \left[ { T_\la - \E T_\la
\over ( \Var T_{ \la} )^{1/2}  } \leq t \right] - \Phi(t) \right|
\leq C ( \log \la )^{dq}  \la ( \Var T_{ \la} )^{-q/2}. \ee
% and
%  \be \label{rates3}   \left| P \left[ {
%T_\la  \over ( \Var T_{ \la} )^{1/2}  } \leq t \right]
%- \Phi(t) \right| \leq  C'(d) ( 1 + |t|)^{-3}
%  M^{3 + \de}
% ( \log \la )^{2d}  \la ( \Var T_{ \la} )^{-3/2} + C'(d)\la^{-2} . \ee
\end{theo}

% The `usual' case of the above theorem, arising in many applications,
%has (\ref{mom}) holding for some $p >3$ so that (\ref{rates2}) holds
%with $q=3$.

Separate arguments are required to establish the asymptotic
behavior of
 the denominator
$(\Var(T_\la))^{1/2}$ in (\ref{rates2}). When $A_\la =A$ for all
$\la$,  it is typically the case for polynomially stabilizing
functionals satisfying moments conditions along the lines of
(\ref{mom}) that there is a constant $\sigma^2 (f,\xi,\kappa) \geq
0$ such that
 \bea
\label{BYeq}
\lim_{\la \to \infty}
\la^{-1}  \Var
\langle f, {\mu}^{\xi}_{\la } \rangle   =
%\int_A f(x)^2 V^\xi(\kappa(x)) \kappa (x) dx
\sigma^2(f,\xi,\ka).
\eea
For further information about $\sigma^2(f,\xi,\ka)$ and
 precise conditions under which (\ref{BYeq}) holds,
see Theorem 2.4(i) of \cite{BY2}.
% with $V^\xi(\cdot)$ given explicitly in terms of $\xi$ in \cite{BY2}.
When (\ref{BYeq}) holds,
by combining it
 with  Theorem \ref{mainthm1} we obtain
\bea
\langle f, \la^{-1/2}\overline{\mu}_{\la }^{\xi}\rangle
%\tod \NN \left(0, \int_A f(x)^2 V^\xi(\kappa(x)) \kappa (x) dx\right),
\tod \NN \left(0, \sigma^2 (f,\xi,\kappa)\right) ,
\label{CLTeq}
\eea
where $\NN(0,\sigma^2)$ denotes a centered normal distribution with
variance $\sigma^2$ if $\sigma^2 >0$, and a unit point mass at 0
if $\sigma^2=0$.

 In many applications (\ref{BYeq}) holds
with $\sigma^2(f,\xi,\ka)
>0$, showing that the case $q=3$ of (\ref{rates2})
 yields
a rate of convergence $O(( \log \la)^{3d} \la^{-1/2} )$ to the
normal distribution.  In other words, we will make frequent use
of:

\begin{coro}\label{maincoro}
Suppose $\|\kappa\|_\infty < \infty$. Suppose that
  $\xi$ is exponentially stabilizing and
satisfies the moments condition (\ref{mom}) for some $p > 3$.  Let
$f \in B(A)$ and put $T_\la := \langle f, \mu^{\xi}_{\la }
\rangle$. If (\ref{BYeq}) holds with
 $\sigma^2 (f,\xi,\kappa) > 0$, then there exists
 a finite constant $C$ depending on $d, \xi$,  $\ka$, $p$ and $f$,
such that for all $\la \geq 2$,
 $$
 \sup_{t \in \R}  \left| P \left[ { T_\la - \E
T_\la \over ( \Var T_{ \la} )^{1/2}  } \leq t \right] - \Phi(t)
\right| \leq C ( \log \la )^{3d}  \la^{-1/2}.
$$
\end{coro}

Our methods actually yield normal approximation and a central
limit theorem when the exponential decay condition is replaced by
a polynomial decay condition of sufficiently high order. We give a
further result along these lines.
\begin{theo}
\label{polyth}
Suppose $\|\kappa\|_\infty < \infty$. Suppose
 for some $p >3$ that
  $\xi$ is polynomially stabilizing of order $\gamma$
with $\gamma > d(150 + 6/p)$, and
satisfies the moments condition (\ref{mom}).
Let $f \in B(A)$ and put $T_\la :=
\langle f, \mu^{\xi}_{\la } \rangle$.
Suppose that
%$\la^{-1} \Var(T_\la)$ converges to a non-zero limit.
(\ref{BYeq}) holds for some $\sigma^2 \geq 0$.
  Then (\ref{CLTeq}) holds and if $\sigma^2 :=\sigma^2(f,\xi,\kappa) >0$
there exists
 a finite constant $C$ depending on $d, \xi$,  $\ka$, $p$ and $f$,
such that for all $\la \geq 2$,
 \be
\label{rates2a} \sup_{t \in \R}  \left| P \left[ { T_\la - \E T_\la
\over ( \Var T_{ \la} )^{1/2}  } \leq t \right] - \Phi(t) \right|
\leq C   \la^{(150pd + 6d - p \gamma )/2(p \gamma - 6d)} . \ee
%In particular, $\la^{-1/2}(T_\la - \E T_\la)$ is asymptotically normal
%provided $\gamma > d( 150 + 6/p)$.
\end{theo}

{\bf Remarks}
\begin{enumerate}
\item
%(i)
 Our results are stated for {\em marked} Poisson point processes,
i.e., for Poisson processes in $\R^d \times \MM$ where $\MM$ is the mark
space. These results are reduced to the corresponding
 results for unmarked Poisson point
processes in $\R^d$ by taking $\MM$ to have a single element (denoted $m$, say)
and identifying $\R^d \times \MM$ with $\R^d $ in the obvious way
by identifying $(x,m)$
with $x$ for each $x \in \R^d$. In this case the notation
(\ref{stab2}) is particularly appropriate. Other
treatments such  as \cite{BY2,Pe2,PY1,PY2}
tend to concentrate on the unmarked case
with  commentary that the proofs carry through to the marked case;
here we spell out the results and proofs in the more general marked case, which
seems worthwhile since examples such as those in Section \ref{seqsubsec}
use the results for marked point processes. Our  examples
in Sections \ref{subsecnng}, \ref{secvoro}, and \ref{subsecoffl}
refer to unmarked point processes and in these examples
we identify  $\R^d \times \{m\}$ with $\R^d$
as indicated above (so that $\Po_\la$ is viewed as a Poisson process in $\R^d$).

\item
 We are not sure
 if the logarithmic factors can be removed in Theorem \ref{mainthm1} or Corollary \ref{maincoro}.
Avram and Bertsimas \cite{AB} obtain a rate of $O((\log
\lambda)^{1+3/(2d)}\lambda^{-1/4})$, for the
length of the $k$-nearest neighbors (directed) graph, the Voronoi
graph, and the Delaunay graph (see Sections \ref{subsecnng}
and \ref{secvoro}). Our method for general stabilizing
functionals is based on theirs, but uses
a stronger general normal approximation result (Lemma \ref{ChenShao} below).

\item
%(iii)
 If (\ref{BYeq}) holds with $\sigma^2(f,\xi,\ka) =0$,
then (\ref{CLTeq}) holds trivially by Chebyshev's inequality, but
Theorem \ref{mainthm1} does not provide any useful information on
rate of convergence. In examples of interest, it can usually be established
that $\sigma^2(f,\xi,\ka) >0$, by further separate arguments.
We do not discuss these in detail here but refer the reader to
\cite{AB,PY1,BY2}.

\item
%(iv)
 Theorems \ref{mainthm1}, \ref{polyth}, and Corollary
\ref{maincoro} require neither the underlying density function
$\kappa$ nor the test function $f$ to be continuous (both of these
conditions are imposed in \cite{BY2}). In particular,
these three results
% Theorems \ref{mainthm1}, \ref{polyth}, and Corollary \ref{maincoro}
 apply when
  $f$ is the indicator function of a Borel subset $B$
of $A$, giving normal approximation for $\bar{\mu}_\la^\xi(B)$.

\item
%(v)
We do not have rate of
convergence results in the binomial (non-Poisson) setting. For
central limit theorems in the binomial setting, we refer to
\cite{PY1} and \cite{BY2}, which treat uniform and non-uniform
samples respectively.

\item
 Some functionals, such as  those defined in terms of the
minimal spanning tree, stabilize without any known bounds on the
rate of decay of the tail probability $\tau(t)$. In these cases
univariate and multivariate central limit theorems hold
 \cite{Pe2,PY1} but our Theorems
\ref{mainthm1} and \ref{polyth} do not apply and explicit rates of
convergence are not known.

\end{enumerate}

\section{Applications}
\label{secapp}

\allco

Applications of Corollary  \ref{maincoro} to
 geometric probability include
functionals of proximity graphs,
 germ-grain models, and random
sequential packing models. The following examples are for
illustrative purposes only and are not meant to be encyclopedic.
For simplicity we will assume that $\R^d$ is equipped with the
usual Euclidean metric. While translation invariance is not needed
in the general results in Section \ref{secgr},
% Corollary \ref{maincoro},
 most of the examples treated in this section
involve translation invariant functionals $\xi$.  However, the
examples can be modified to treat the (non-translation-invariant)
situation where $\R^d$ has  a  local metric structure.

\subsection{$k$-nearest neighbors graph}
\label{subsecnng}

 Let $k$ be  a positive
integer. Given a locally finite  point set $\X \subset \R^d$,
  the $k$-nearest neighbors (undirected) graph
  on $\X$, denoted $kNG(\X)$, is the graph with vertex set
$\X$ obtained by including $\{x,y\}$ as an edge whenever $y$ is
one of the $k$ nearest neighbors of $x$ and/or $x$ is one of the
$k$ nearest neighbors of $y$. The $k$-nearest neighbors (directed)
graph on $\X$, denoted $kNG'(\X)$, is the graph with vertex set
$\X$ obtained by placing a directed edge between each point and
its $k$ nearest neighbors.

Let $N^k(\X)$ denote the total edge length of the (undirected)
$k$-nearest neighbors graph on $\X$.  Note that $N^k(\X) = \sum_{x
\in \X} \xi^k(x; \X)$, where $\xi^k(x; \X)$ denotes half the sum of the
edge lengths in $kNG(\X)$ incident to $x$.
  If $A$ is convex or polyhedral and $\ka$ is bounded away from $0$
on $A$, then $\xi^k$ is exponentially
stabilizing (cf. Lemma 6.1 of \cite{PY1})
%the add one cost is non-degenerate (Lemma 6.3 of \cite{PY1}),
 and has moments of all orders. Moreover,
as shown in \cite{BY2} (see e.g. display (2.11), Theorem 3.1),
at least when $f$ is continuous and $A_\la =A$ for all $\la$,
 $$
 \lim_{\la \to \infty} \la^{-1}  \Var \langle f,
{\mu}^{\xi}_{\la } \rangle   = V^\xi \int_A f(x)^2
\kappa(x)^{(d-2)/d} dx,
$$
where $V^\xi$ denotes the limiting variance for the total edge
length of the $k$-nearest neighbors graph on $\la^{1/d} \P_\la$
when $\ka$ is the uniform distribution on the unit cube.  Since
%$V^\xi$ is strictly positive (Lemma 6.3 of \cite{PY1}), it follows
$V^\xi$ is strictly positive (Theorem 6.1 of \cite{PY1}), it follows
that (\ref{BYeq}) holds with $\sigma^2(f,\xi^k, \kappa) >0$.
 We thus obtain via Corollary \ref{maincoro} the following rates
in the CLT
for the total edge length of $N^k(\la^{1/d} \P_\la)$
%, in the case where $A$ is a cube,
improving  upon Avram and Bertsimas \cite{AB} and Bickel and
Breiman \cite{BB}.  A similar CLT holds for the total edge length
of the $k$-nearest neighbors directed graph.

\begin{theo}
Suppose $A$ is convex or polyhedral and $\ka$ is bounded away
from $0$ on $A$.
 Let $N_\la := N^k(\la^{1/d} \P_\la)$ denote the total edge
length of the $k$-nearest neighbors graph on $\la^{1/d} \P_\la$.
There exists a finite constant $C$ depending on $d, \xi^k,$ and
$\ka$ such that
 \be
\label{nnrates2} \sup_{t \in \R}  \left| P \left[ { N_\la - \E
N_\la \over ( \Var N_{ \la} )^{1/2}  } \leq t \right] - \Phi(t)
\right| \leq   C ( \log \la )^{3d}  \la^{-1/2}. \ee
% and
%  \be \label{nnrates3}   \left| P \left[ {
%N_\la - \E N_\la \over ( \Var N_{ \la} )^{1/2}  } \leq t \right] -
%\Phi(t) \right| \leq  C'(d) ( 1 + |t|)^{-3}
%    ( \log \la )^{2d}  \la^{-1/2}. \ee
\end{theo}

Similarly, letting
$\xi^s(x; \X)$ be one or zero  according to whether the distance
between $x$ and its nearest neighbor in $\X$ is less than $s$ or
not, we can verify that $\xi^s$ is exponentially stabilizing and that
%the add-one cost is non-degenerate,
the variance of $\sum_{   x \in \la^{1/d} \P_\la }
 \xi^s(x; \la^{1/d} \P_\la )$ is bounded below by a positive multiple of
$\la$.  We thus obtain rates of convergence
of $O ( ( \log \la )^{3d}  \la^{-1/2} ) $
in the CLT for the one-dimensional  marginals of
 the empirical distribution function of $k$ nearest neighbors distances on
$\la^{1/d} \P_\la$,
improving upon those implicit on p. 88 of \cite{Pe}.

Using the results from section 6.2 of
\cite{PY1}, we could likewise obtain the same rates of convergence in the CLT for
%the total number of components in the $k$ nearest neighbors graph as well
the number of vertices of fixed degree in the $k$ nearest neighbors graph.

Finally in this section, we re-consider the non-translation-invariant
%\Comment{new para. MP}
 example given in Section \ref{secgr}, where a point at $x$
is colored red with probability $q(x)$ and green with probability
$1-q(x)$, and $\xi(x;\X)$ takes the value 0 if the nearest neighbor
of $x$ in $\X$ has the same color as $x$, and takes the value 1 otherwise.
We can use Corollary \ref{maincoro} to derive a central limit theorem,
with $O((\log \la)^{3d}\la^{-1/2})$ rate of convergence,
for $\sum_{x \in \Po_\la} f(x) \xi(x;\Po_\la)$, where
$f$ is a bounded measurable test function.

\subsection{Voronoi
%, Delaunay,
 and sphere of influence graphs}
\label{secvoro}

We will consider the Voronoi graph
for $d =2$ and the
% and its planar
%dual, the Delaunay graph; for $d \geq 3$ we will consider the
sphere of influence graph for all $d$. Given a locally finite set
$\X \subset \R^2$ and given $x \in \X$, the locus of points closer
to $x$ than to any other point in $\X$ is called the {\em Voronoi}
cell centered at $x$. The graph consisting of all boundaries of
Voronoi cells is called the Voronoi graph generated by $\X$.

The sum of the lengths of the finite edges  of the Voronoi graph
on $\X$ admits the representation $\sum_{x \in \X} \xi(x; \X)$,
where $\xi(x; \X)$  denotes one half the sum of the lengths of the
finite edges in the Voronoi cell at $x$.  If  $\ka$ is bounded
away from $0$ and infinity and $A$ is convex, then geometric
arguments show that there is a random variable $R$ with
exponentially decaying tails such that for any $x \in \cal P_\la$,
the value of $\xi(x; \P_\la)$ is unaffected by points outside
$B_{\la^{-1/d} R}(x)$ \cite{BY2,PY1,PY4}.
%\Comment{Was $B_R$. MP}
 In other words,  $\xi$ is
exponentially stabilizing and satisfies the moments condition
(\ref{mom}) for all $p > 1$.
 Also, the variance
of the total edge length of these graphs on $\la^{1/d}\P_\la$
%\Comment{Was $\P_\la$. MP}
 is bounded
below by a multiple of $\la$. We thus obtain
 $O ( ( \log \la )^{3d}  \la^{-1/2} ) $
rates of convergence
 in the CLT for the
total edge length functionals of these graphs on $\la^{1/d}\P_\la$,
%\Comment{Was $\P_\la$. MP}
thereby
improving and generalizing the results of Avram and Bertsimas
\cite{AB}.

Given a locally finite set  $\X \subset \R^d$,
 the sphere of influence graph $\SIG(\X)$ is
a graph with vertex set $\X$, constructed as follows:  for each $x
\in \X$ let $B_x$ be a ball around $x$ with radius equal to
$\min_{y \in \X \setminus \{x\}}
 \{\vert y - x \vert \}.$  Then $B_x$ is
called the {\em sphere of influence} of $x$.  We put an edge
between $x$ and $y$ iff the balls $B_x$ and $B_y$ overlap. The
collection of such edges is the {\em sphere of influence graph}
(SIG) on $\X$.

The total number of edges of the sphere of influence graph on $\X$
admits the representation $\sum_{x \in \X} \xi(x; \X)$, where
$\xi(x; \X)$ denotes one half the degree of SIG at the vertex $x$.
The number of vertices of fixed degree $\de$ admits a similar
representation, with $\xi(x; \X)$ now equal to one (respectively,
zero) if the degree at $x$ is $\de$ (respectively, if degree at
$x$ is not $\de$). If  $\ka$ is bounded away from $0$ and infinity
and $A$ is convex, then geometric arguments show that both choices
of the functional $\xi$ stabilize (see sections 7.1 and 7.3 of
\cite{PY1}).
 Also, the variance
of both the total number edges and the number of vertices of fixed
degree in the SIG on $\la^{1/d} \P_\la$
%\Comment{Was $\Po_\la$}
 is bounded below by a multiple of
$\la$ (sections 7.1 and 7.3 of \cite{PY1}). We thus
 obtain  $O ( ( \log \la )^{3d}  \la^{-1/2} ) $
rates of convergence
 in the CLT for the total number of edges and the
number of vertices of fixed degree in the
sphere of influence graph on $\P_\la$.

\subsection{Random sequential packing models}
\label{seqsubsec}

 The following
prototypical random sequential packing model is of considerable
scientific interest; see \cite{PY2} for references to the vast
literature.

With $N(\la)$ standing  for a Poisson \ random variable with
parameter $\la$, we  let $B_{\la,1},B_{\la,2},...,B_{\la,N(\la)}$
be
 a sequence
 of $d$-dimensional balls of volume $\la^{-1}$  whose centers
are  i.i.d. random $d$-vectors $X_1,...,X_{N(\la)}$ with probability
density function $\ka: A \to [0,\infty) $.
 Without loss
of generality, assume that the balls are sequenced in the order
determined by marks (time coordinates) in $[0,1]$.  Let the first
ball $B_{\la,1}$ be {\em packed}, and recursively for $i=2,3,
\ldots, $ let
 the $i$-th ball
$B_{\la,i}$ be  packed iff $B_{\la,i}$ does not overlap any ball
in $B_{\la,1},...,B_{\la,i-1}$ which has already been packed. If
not packed, the $i$-th ball is discarded.

Packing models of this type arise in diverse disciplines,
including physical, chemical, and biological processes \cite{PY2}.
Central limit theorems for the number of accepted (i.e., packed) balls are
established in \cite{PY2, BY2}, whereas laws of large numbers are
given in \cite{PY4}.

Let $\MM=[0,1]$ with $\PP_\MM$ being the uniform distribution on the
unit interval. For any finite point set  $\XX \subset \R^d \times
[0,1]$,
 assume the points $(x,s)
\in \XX$  represent the locations and arrival times.
  Assume balls of
volume $\la^{-1}$ centered at the locations of $\XX$ arrive
sequentially in an order determined by the time coordinates, and
assume as before that each ball is packed or discarded according
to whether or not it overlaps a previously packed ball.  Let
$\xi((x,s); \XX)$ be either $1$ or $0$ depending on whether the ball
centered at $x$ at times $s$ is packed or discarded. Consider the re-scaled
packing functional $\xi_\la((x,s);\X) = \xi((\la^{1/d} x,s); \la^{1/d}
\X)$,
where balls centered at
points of $\la^{1/d} \X$ have volume one.  The random measure
$$
\mu^{\xi}_{\la} := \sum_{i=1}^{N(\la)} \xi_{\la} ((X_i,U_i); \{(X_j,
U_j)
\}_{j=1}^{N(\la)}  ) \ \de_{X_i},
$$
is called the random sequential packing measure induced by balls
with centers arising from $\ka$. The convergence of the finite
dimensional distributions of the packing measures
$\mu^{\xi}_{\la}$ is established in \cite{BY1,BY2}.
  $\xi$ is exponentially
stabilizing \cite{PY2,BY1} and for any continuous $f \in B([0,1]^d)$ and
$\ka$ uniform, the variance of $\langle f, \mu^{\xi}_{\la}\rangle $
 is bounded below by a positive multiple of $\la$ \cite{BY2},
showing that $\langle f, \mu^{\xi}_{\la}\rangle $ satisfies a CLT
with an $O ( ( \log \la )^{3d}  \la^{-1/2} ) $ rate of
convergence.

%MP. Do the BY papers establish the results claimed in the prev.
%paragraph only for $\mu(f)$ with $f$ continuous, and hence not
%necessarily for $\mu(B)$? If so, adapt the preceding para.
%accordingly. (This is now changed.  JY)

It follows easily from the stabilization analysis of  \cite{PY2}
that many variants of the above basic packing model satisfy
similar rates of convergence in the CLT.
 Examples include balls of bounded random radius,  cooperative sequential adsorption
 (\cite{PY2}), and
monolayer ballistic deposition (\cite{PY2}). In each case the
number of particles accepted satisfies the CLT with an $O ( ( \log
\la )^{3d} \la^{-1/2} ) $ rate of convergence. The same comment
applies for the number of seeds accepted in spatial birth-growth
models \cite{PY2}.

\subsection{Independence number, off-line packing}
\label{subsecoffl}

An {\em independent set} of vertices in a graph $G$ is
a set of  vertices in $G$, no two of which
are connected by an edge.
The {\em independence number} of $G$, which we denote $\beta(G)$,
is defined to be the maximum cardinality of all
 independent sets of vertices in $G$.

For $r>0$, and for finite or countable $\X \subset \R^d$,
let $G(\X,r)$ denote the
{\em geometric graph} with vertex set $\X$ and with edges between each pair
 of vertices distant at most  $r$ apart.
Then the independence number $\beta(G(\X,r))$ is the
maximum number of disjoint closed balls of radius $r/2$ that can be
centered at points of $\X$; it is an `off-line' version of
the packing functionals considered in the previous section.

Let $b>0$ be a constant, and
 consider the graph $G(\Po_\lambda,b\lambda^{-1/d})$
(or equivalently, $G(\lambda^{1/d} \Po_\lambda,b)$).
Random geometric graphs of this type are the subject
of \cite{Pe}, although independence number is
considered only briefly there (on page 135). A law of large
numbers for the independence number is described in
Theorem  2.7 (iv) of \cite{PY4}.

  For $u>0$, let $\H_u$ denote a homogeneous Poisson point process
of intensity $u$ on $\R^d$, and let $\H_{u}^0$ be the point
process $\H_u$
 with a point inserted at the origin.
As on page 189 of \cite{Pe},
let $\lambda_\rc$ be the infimum of all
$u$ such that the origin has a non-zero probability
of being in an infinite component of $G(\H_u,1)$.

If $b^d \| \kappa\|_\infty < \lambda_{\rc}$, we can use Corollary
 \ref{maincoro} to
 obtain a central limit theorem
for the independence number
$\beta(G(\Po_\lambda,b \lambda^{-1/d}))$,
namely
\bea
\label{offclt}
 \sup_{t \in \R}  \left| P \left[ {
\beta(G(\Po_\lambda,b \lambda^{-1/d})) -
 \E
 \beta(G(\Po_\lambda,b \lambda^{-1/d}))
 \over ( \Var
 \beta(G(\Po_\lambda,b \lambda^{-1/d}))
 )^{1/2}  } \leq t \right] - \Phi(t)
\right| \leq C ( \log \la )^{3d}  \la^{-1/2}.
\eea

We sketch  the proof. For finite $\X \subset \R^d$ and $x \in \X$,
  let $\xi(x;\X)$ denote the
independence number of the component of $G(\X,b)$ containing
vertex $x$,  divided by the number of vertices in this component.
Then $\sum_{x \in \X} \xi(x;\X)$ is the independence number of
$G(\X,b)$, since the independence number of any graph is the sum
of the independence numbers of its components. Also, our choice of
$\xi$ is translation-invariant, and so we obtain \bean \beta
(G(\Po_\la, \la^{-1/d}b)) = \beta( G(\la^{1/d} \Po_\la, b) ) =
\sum_{x \in \Po_\la} \xi(\la^{1/d} x;\la^{1/d} \Po_\la)
\\
= \sum_{x \in \Po_\la} \xi_\la(x; \Po_\la) = \langle \mu_\la^\xi,
f \rangle, \eean where we here take the test function $f$ to be
identically 1 and take $A_\la =A$ for all $\la$. Thus a central
limit theorem holds for $\beta (G(\Po_\la, \la^{-1/d}b))$ by
application of Corollary
 \ref{maincoro}, if  $\xi$ and $\kappa$ satisfy the conditions
for that result.

We take $R(x,\la)$ to be the distance from $\la^{1/d}x$ to the
furthest point in the component containing $\la^{1/d}x$ of
$G(\la^{1/d}\Po_\la, b)$, plus $2b$.
 Since $\xi_\la(x; \Po_\la)$ is
determined by
% the independence number of
the component of $G(\la^{1/d}\Po_\la,b)$
  containing $\la^{1/d} x$,
%divided by the number of elements of this component,
 and this component  is unaffected by the addition
or removal of points to/from $\Po_\la$ at a distance greater than
$\la^{-1/d}R(x,\la)$ from $x$, it is indeed the case that
  $R(x,\la)$ is a radius of stabilization.

The point process
 $\lambda^{1/d}\Po_\lambda$
is dominated by $\H_{\|\kappa\|_\infty}$ (in the sense of
\cite{Pe}, page 189).
Hence, $P[R(x,\la) > t]$ is bounded by the probability
that the component containing $x$ of
 $G(\H_{\|\kappa\|_\infty}\cup \{ \la^{1/d}x \},b)$
has at least one vertex
outside $B_{t-2b}(\la^{1/d}x)$.  This probability does not depend on $x$, and equals
the probability that
the component of $G(\H^0_{b^d \|\kappa\|_\infty},1)$ containing
the origin includes a vertex outside $B_{(t/b)-2}$.
By exponential decay
for subcritical continuum percolation (Lemma 10.2 of \cite{Pe}),
this probability decays exponentially in $t$, and exponential
 stabilization of $\xi$ follows. The moments condition (\ref{mom})
is trivial  in this case, for any $p$, since $0 \leq \xi(x;\X) \leq 1$.

 Thus, Corollary
 \ref{maincoro} is indeed applicable, provided that (\ref{BYeq})
holds in this case, with $\sigma^2 >0$. Essentially (\ref{BYeq})
follows from Theorem 2.1 of \cite{BY2}, with strict inequality
$\sigma^2 >0$ following from (2.10) of \cite{BY2};
  in the case where $\kappa$ is the density
function of a uniform distribution on some suitable subset of
$\R^d$, one can alternatively use Theorem 2.4 of \cite{PY1}. We do
not go into details here about the application of results in this
example, but we do comment further  on why the distribution of the
`add one cost' (see \cite{PY1,BY2}) of insertion of a point at the
origin into a homogeneous Poisson process is nondegenerate, since
this is needed to verify $\sigma^2
>0$ and this  example was not considered in \cite{PY1} or
\cite{BY2}.

The above add one cost is
 the variable denoted $\Delta(\infty)$
in the notation of \cite{PY1}, or $\Delta^\xi(u)$
%\Comment{changed to $\Delta^\xi(u)$. JY}
  in the notation of \cite{BY2}. It
%The variable $\Delta(\infty)$
is the independence number
of the component containing the origin of $G(\H_u;b)$ minus
the independence number of this component with the origin removed
(we need only to consider the case where $b^du$ is subcritical).
This variable can take the value 1, for example if
the origin is isolated in $G(\H_u;b)$, or zero, for example if
 the component containing the origin has two vertices.
Both of these possibilities  have  strictly positive probability,
and therefore  $\Delta(\infty)$ has a non-degenerate distribution.

\allco

\section{Proof of Theorems}

\allco

\subsection{A CLT for dependency graphs}

We shall prove Theorem \ref{mainthm1} by showing that exponential
stabilization implies  that a modification of
 $\langle f, \overline{\mu}^{\xi}_{\la } \rangle$ has a
 {\em dependency graph} structure,
 whose definition we now recall  (see e.g. Chapter 2 of
\cite{Pe}).  Let $X_{\alpha}, \ \alpha \in {\cal V},$ be a
collection of random variables. The graph $ G := ({\cal V}, {\cal
E } )$ is a {\em dependency graph} for $X_{\alpha}, \ \alpha \in
{\cal V},$ if for any pair of disjoint sets $A_1, A_2 \subset
{\cal V}$ such that no edge in ${\cal E }$ has one endpoint in
$A_1$ and the other in $A_2$, the sigma-fields $\sigma\{
X_{\alpha}, \alpha \in A_1 \}$, and $\sigma\{ X_{\alpha}, \alpha
\in A_2 \}$, are mutually independent. Let $D$ denote the maximal
degree of the dependency graph.

It is well known that sums of random variables indexed by the
vertices of a dependency graph admit rates of convergence to a
normal.  The rates of Baldi and Rinott \cite{BR} and those in
Penrose \cite{Pe} are particularly useful; Avram and Bertsimas
\cite{AB} use the former to obtain rate results for the total edge
length of the nearest neighbor, Voronoi, and Delaunay graphs.

In many cases, the following theorem of Chen and Shao \cite{CS}
provides superior rate results.
%  We shall apply this result when $p = 3$.
%\Comment{removed sntnce. MP}
 For any random variable $X$ and any $p > 0$, let
 $||X||_p = (\E[|X|^p ])^{1/p}.$

\begin{lemm} \label{ChenShao}  (see Thm 2.7 of \cite{CS})
Let $2 < q \leq 3$.  Let $W_i, \ i \in {\cal V},$ be random
variables indexed by the vertices of a dependency graph. Let $W =
\sum_{i \in {\cal V} } W_{i}.$ Assume  that $\E [W^2] = 1, \E
[W_i] = 0,$ and $||W_i ||_q \leq \theta$ for all $i \in  {\cal V}$
and for some $\theta > 0.$ Then \be
 \label{CS1} \sup_t | P [ W \leq t ] - \Phi(t) | \leq 75 D^{5(q
- 1) } |  {\cal V} | \theta^q.
 \ee
\end{lemm}

\vskip.5cm

\subsection{Auxiliary lemmas}
%\Comment{maybe need $A$ compact convex w. non-empty interior?}
 To
prepare for the proof of Theorem \ref{mainthm1} we will need some auxiliary
lemmas.
 Throughout, $C$  denotes a generic
constant depending possibly on $d, \ \xi,$ and $\kappa$ and whose
value may vary at each occurrence. We assume $\la > 1$ throughout.

Let $(\rho_\la, \la >0)$ be a function to be chosen later,
in such a way that $\rho_\la \to \infty$ and $\la^{-1/d} \rho_\la \to 0$
as $\la \to \infty$.  Given $\la>0$, let  $s_\la:= \lambda^{-1/d} \rho_\la$,
and let $V := V(\la)$ denote the number of cubes of the form
$Q= \prod_{i=1}^d [j_i s_\la, (j_i+1)s_\la)$,
 with all $j_i \in \Z$, such that $\int_{Q}\ka(x)dx >0$;
enumerate these cubes as $Q_1,Q_2,\ldots,Q_{V(\la)}$.
Since $\ka$ is assumed to have bounded support,
 it is easy to see that $V(\la) = O(\la \rho_\la^{-d})$
as $\la \to \infty$.

For all $1 \leq i \leq V(\la)$, the number of points of $\P_\la
\cap (Q_i\times \MM)$ is a Poisson random variable $N_i:= N(\stau_i)$, where
\bea
\stau_i := \la \int_{Q_i} \ka(x)dx  \leq \| \ka \|_\infty \rho_\la^d.
\label{nudef}
\eea
 Assuming $\stau_i > 0$, choose
an ordering on the points of $\P_\la \cap (Q_i\times \MM)$
 uniformly at random
from all $(N_i)!$ possible such orderings. Use this ordering to
list the points as $(X_{i,1},U_{i,1}),...,(X_{i,N_i},U_{i,N_i})$,
 where conditional on
the value of $N_i$, the random variables
 $X_{i,j}, \ j = 1,2,...$ are i.i.d.
on $Q_i$ with a density  $\ka_i(\cdot) := \ka( \cdot)/ \int_{Q_i}
\ka(x) dx$, and the $U_{i,j}$ are i.i.d. in $\MM$ with
distribution $\PP_\MM$, independent of $\{X_{i,j}\}$. Thus we have
the representation
  $\P_\la =  \cup_{i=1}^{V(\la)} \{(X_{i,j},U_{i,j}) \}_{j=1}^{N_i} $.
  For all $1 \leq i \leq
V(\la)$, let $\P_i := \P_\la \setminus \{(X_{i,j},U_{i,j})
\}_{j=1}^{N_i}$ and note that $\P_i$ is a Poisson point process on
$\R^d \times \MM$ with intensity
 $\la \ka(x) {\bf 1}_{\R^d \setminus Q_i}(x) dx \times \PP_{\MM} (ds) $.
 %on $[0, 1]^d \setminus Q_i$ and
%intensity zero on $Q_i$.

We show that the condition (\ref{mom}),
which
bounds the moments of the value of $\xi$ at points inserted into $\P_\la$,

 also yields bounds on
$\E [ | \xi_{\la}((X_{i,j},U_{i,j});  \P_{\la} )|^p
%\cdot
\1_{A_\la}(X_{i,j})
\1_{j \leq N_i} ]$. More precisely, we have:

\begin{lemm} \label{lemMbd}
 If (\ref{mom}) holds for some $p>0$,
%$M^p < \infty$, $p > 0,$
 then there is a
constant $C:=C(p)$ %depending on $M^p$ and $\ka$
 such that for all $\la > 1$, all $j \geq 1$ and
 $1 \leq i \leq V(\la)$,
 \be \label{Mbound}
 \E [ | \xi_{\la}(X_{i,j};  \P_{\la} )
 \cdot
{\bf 1}_{A_\la}(X_{i,j})
\1_{j \leq N_i}  |^p ] \leq C \rho_\la^d. \ee
\end{lemm}
{\em Proof.}
 If $N_i = n$, then denote $\{(X_{i,1},U_{i,1}),...,(X_{i,N_i},U_{i,N_i})\}$ by
$\X_n$. By definition,
$$
\E [ | \xi_{\la}((X_{i,j},U_{i,j});  \P_{\la} )
 \cdot
{\bf 1}_{A_\la}(X_{i,j})
 \1_{j \leq N_i} |^p ] =
\sum_{n = j}^{\infty}  \int_{Q_i \cap A_\la}
\E [ | \xi_\la((x,U); \X_{n-1} \cup
\P_i) |^p ] \ka_i(x) dx
  \cdot P[N_i = n],
$$
where the expectation on the right hand side is
with respect to $U$, $\X_{n-1}$ and $\P_i$. The above is bounded by

$$
\leq \stau_i \sum_{n =1}^{\infty} \int_{Q_i\cap A_\la}
\E [| \xi_\la((x,U); \X_{n-1}
\cup \P_i) |^p ] \ka_i(x) dx   \cdot {e^{- \stau_i} \stau_i^{n-1}
\over (n-1)! }
$$
$$
= \stau_i \sum_{m = 0}^{\infty}  \int_{Q_i\cap A_\la} \E \left[ | \xi_\la((x,U);
\P_{\la}) |^p  \  \ |  \ \ | \P_\la \cap (Q_i \times \MM)| = m \right]
\ka_i(x) dx  \cdot P[| \P_\la \cap ( Q_i \times \MM)| = m]
$$
$$
= \stau_i  \int_{Q_i\cap A_\la} \E [ | \xi_\la((x,U);  \P_{\la}  ) |^p ]
\ka_i(x) dx \leq
 {\rm const.} \times \stau_i,
 %\stau_i M^p,
$$
 where the last inequality follows
by (\ref{mom}).
%Since $\stau_i = \la \int_{Q_i} \ka(x)dx
%\leq \|\ka\|_{\infty} \rho_\la^d$,
 By (\ref{nudef}), this shows
(\ref{Mbound}). \qed

\ \

%For all $1 \leq i \leq V$ and for $\stau_i := \la \int_{Q_i}
%\ka(x)dx$, let $N_i := N(\stau_i)$ be an independent Poisson random
%variable with parameter $\stau_i$.
For $1 \leq i \leq V$, and for  $j \in \{ 1,2,\ldots \}$, we define

  $$
%\xi_j :=
\xi_{i,j} := \left\{ \begin{array}{ll}
\xi_\la ((X_{i,j},U_{i,j});  \P_\la )
& {\rm if }~  N_i \geq j, X_{i,j} \in A_\la \\
0 & {\rm otherwise}
\end{array}
\right.
 %\cdot 1_{ E_{i,j} }
$$
% when $1 \leq j \leq N_i$ and
%otherwise we set $\xi_j = 0$.

\begin{lemm} \label{sumbound}
If (\ref{mom}) holds for some $p>q >1$, then
there exists $C := C(p,q)$ such that
for $1 \leq i \leq V(\la)$,
\bea
\left| \left|  \sum_{j = 1}^{\infty} |\xi_{i,j} | \right| \right|_q
\leq C  \rho_\la^{d(p+1)/p}.
\label{bound3}
\eea
\end{lemm}
{\em Proof.}
Fix $i \leq V(\la)$ and write $\xi_j $ for $\xi_{i,j}$.
Clearly, with $N := N_i$ and $\stau := \stau_i$,
$$
\left| \left| \sum_{j =1}^{\infty} |\xi_j | \right| \right|_q =
\left| \left| \sum_{j=1}^{\infty}| \xi_j|
 \left(  \1_{ N
\leq \stau }  + \sum_{k = 0}^{\infty}  \1_{2^k \stau < N \leq 2^{k +
1} \stau } \right)
 \right| \right|_q
$$
$$
 \leq \left| \left| \sum_{j=1}^{\infty}   \sum_{k = 0}^{\infty} |
\xi_j | \cdot  \1_{2^k \stau < N \leq 2^{k + 1} \stau } \right|
\right|_q
 +  \left| \left| \sum_{j =1}^{\infty} | \xi_j | \cdot \1_{ N \leq \stau } \right|
 \right|_q.
$$
Since a.s. only finitely many summands in the double sum are
non-zero, by subadditivity of the norm, the above is bounded by
$$
\leq \sum_{k = 0}^{\infty}  \left| \left| \sum_{j =1}^{\infty} |
\xi_j | \cdot  \1_{2^k \stau < N \leq 2^{k + 1} \stau } \right|
\right|_q + \left| \left| \sum_{j =1}^{\lfloor \stau \rfloor } |
\xi_j| \cdot \1_{ N \leq \stau } \right| \right|_q
$$
$$
\leq \sum_{k = 0}^{\infty} \left|  \left| \sum_{j =1}^{\lfloor
2^{k+1} \stau \rfloor } |  \xi_j | \cdot \1_{ N \geq 2^k \stau}
\right| \right|_q +
 \left|  \left| \sum_{j =1}^{\lfloor \stau \rfloor }  |\xi_j | \cdot \1_{N \leq \stau} \right| \right|_q
 $$
\be
\label{bound}
 \leq \sum_{k = 0}^{\infty}  \sum_{j =1}^{\lfloor 2^{k+1}
\stau \rfloor } \left|  \left|  \xi_j  \cdot \1_{ N \geq 2^k
\stau} \right| \right|_q +
  \sum_{j =1}^{ \lfloor \stau \rfloor }  \left|  \left| \xi_j  \cdot \1_{N \leq \stau} \right|
  \right|_q,
 \ee
 where here and elsewhere  $\lfloor x \rfloor$ denotes the
 greatest integer less than or equal to $x$.

With $\eta := (1/q)-(1/p)$,
H\"older's inequality followed by (\ref{Mbound}) yields
 \be
 \label{bound1}
 || \xi_j \cdot \1_{
N \geq 2^k
 \stau}||_q \leq || \xi_j||_{p } \cdot
(P [N \geq 2^k \stau ])^{\eta }
   \leq C  \rho_\la^{d/p}
(P [N \geq 2^k \stau ])^{\eta }.
 \ee

 Substituting into (\ref{bound}) we obtain
\be \label{normbound}
 \left| \left| \sum_{j =1}^{\infty} |\xi_j|
\right| \right|_q \leq C  \rho_\la^{d/p}  \sum_{k =
0}^{\infty} \stau 2^{k+1} \cdot (P [N \geq 2^k \stau ]
 )^{\eta } +
\sum_{j =1}^{\lfloor \stau \rfloor}  \left|  \left| \xi_j  \cdot
\1_{N \leq \stau} \right|
 \right|_p.
 \ee
By tail estimates for the Poisson distribution
 (see e.g. (1.12) in \cite{Pe}), since $e^2 < 8$ we have
% for $2^k \geq e^2$ (i.e., for $k \geq 3$)
$$
(P[ N \geq 2^k \stau ])^\eta
 \leq \left(\exp \left( - 2^{k-1} \nu \log(2^k) \right) \right)^\eta
= \exp(- k2^{k-1} \eta \nu \log 2 ),
~~~ k \geq 3.
$$
Hence,
$$
(P[ N \geq 2^k \stau ])^{\eta} \leq 2^{-2k}, ~~~ k
\geq \max(3, 2 - \log_2 (\eta \nu)).
$$
Hence, since $\eta$ is a constant, for $\stau > 0$ we have
$$
 \sum_{k = 0}^{\infty} 2^{k+1}   \cdot (P [N \geq 2^k \stau ]
 )^{\eta }
 \leq
\sum_{k < \max(3, 2 - \log_2(\eta \stau))} 2^{k+1}
+ \sum_{k \geq \max(3, 2 - \log_2 (\eta \stau))} 2^{1-k}
$$
$$
 \leq  C \max(1/\stau,1) + 2
$$
 where $C$ does not depend on $\stau$.
 Thus by (\ref{nudef}), the  first term in the
right hand side of (\ref{normbound}) is  $O(  \rho_\la^{d+d/p} )$.
Also, the second term in the right hand side of (\ref{normbound})
is $O(\rho_\la^{d+d/p})$ by Lemma \ref{lemMbd}.
Hence, (\ref{normbound}) implies (\ref{bound3}). \qed

\subsection{Proof of Theorems \ref{mainthm1} and \ref{polyth}}
We prove Theorems \ref{mainthm1} and \ref{polyth} together.
When proving Theorem \ref{mainthm1}
we assume that
 $\xi$ is exponentially stabilizing and
(\ref{mom}) holds for some $p>2$, and we choose
$q\in (2,3]$ with $q<p$.
When
 proving Theorem \ref{polyth} we assume that
$\xi$ is polynomially stabilizing
of order $\gamma $ with $\gamma >d(150+6/p)$, and that
 (\ref{mom})
holds for some $p >3$,
 and we set $q=3$.

Throughout
 this section, we fix $f \in B(A)$ and
set $T_\la := \langle f, {\mu}^{\xi}_{\la } \rangle$.
We follow the setup of the preceding section, with
%$[0,1)^d$
the support of $\kappa$ covered by cubes of side $\la^{-1/d}\rho_\la$,
and we now  choose $\rho_\la$. With the tail probability $\tau(t)$
defined at (\ref{taudef}), we choose $\rho_\la $
%= \alpha \log \la$ with $\alpha$ chosen
in such a way that there is a constant $C$ such that for all $\la
\geq 1$, \bea \rho_\la^{d/p} (\la \tau(\rho_\la))^{(q-2)/(2q)} < C
\la^{-4} ~~~~{\rm and } ~~~ \tau(\rho_\la) < C \la^{-3}
\label{rhodef} \eea and also \bea \rho_\la^d < C \la^{p/(p+2)}.
\label{0513a} \eea In the exponentially stabilizing case (Theorem
\ref{mainthm1}) we achieve this by taking $\rho_\la = \alpha \log
\la$ for some suitably large constant $\alpha$. In the
polynomially  stabilizing case of order $\gamma$ (Theorem
\ref{polyth}), we take $\rho_\la =  C \la^{a}$ with \bea a =
\frac{25 p}{p \gamma - 6d} {\rm ~~~ so ~~~} a \left(
\frac{\gamma}6 - \frac{d}{p} \right) = {25\over 6}. \label{adef}
\eea which implies that (\ref{rhodef}) holds with $q=3$, and that
(\ref{0513a}) holds (to obtain the last conclusion we use our
assumption on $\gamma$, which implies $\gamma > d(25 + 56/p)$.)

For all $1 \leq i \leq V$ and all $j =1,2,...,$ let $R_{i,j}$
denote the radius of stabilization of $\xi$ at $(X_{i,j},U_{i,j})$
 if $1 \leq j \leq N_i$ and $X_{i,j}\in A_\la$;
let $R_{i,j}$ be zero otherwise.

%Let  $E_{i,j} := \{ R_{i,j}  \leq \al \log \la  \}$. Then by
Let  $E_{i,j} := \{ R_{i,j}  \leq \rho_\la  \}$.
Let $E_\la := \cap_{i=1}^{V} \cap_{j=1}^{\infty} E_{i,j}$
Then by Markov's inequality and
standard Palm theory (e.g. Theorem 1.6 in \cite{Pe})

\bea
P[E_\la^c] \leq
\E \left[ \sum_{i=1}^{V} \sum_{j=1}^{N_i} \1_{E_{i,j}^c} \right]
%\nonumber
%\\
 =
\la \int_{A_\la} P[ R(x,\la)
> \rho_\la] \ka( x) dx
\leq \la \tau(\rho_\la).
%\leq C \la^{-33}
\label{0505}
\eea

For each $\la$, and $x \in \R^d$,
 set $f_\la(x):= f(x)\1_{A_\la}(x)$.
Recalling the representation
 $\P_\la = \cup_{i=1}^{V(\la)} \{X_{i,j} \}_{j=1}^{N_i} $,
 we have
$$
T_\la =   \sum_{i=1}^{V(\la)} \sum_{j = 1 }^{N_i} \xi_\la
((X_{i,j},U_{i,j}); \P_\la )  \cdot f_\la(X_{i,j}).
$$
To obtain rates of normal approximation for $T_\la$, it will be be
convenient to consider a closely related sum enjoying more
independence between terms, namely

$$
T'_{\la} := \sum_{i=1}^{V(\la)}  \sum_{j = 1 }^{N_i} \xi_\la
((X_{i,j},U_{i,j});
 \P_\la  ) \cdot \1_{E_{i,j}}   \cdot f_\la(X_{i,j}).
$$
For all $1 \leq i\leq V(\la)$ define
$$
S_i := S_{Q_i}:= (\Var T'_{\la})^{-1/2}   \sum_{j =1}^{N_i}
\xi_\la ((X_{i,j},U_{i,j}); \P_\la )  \cdot \1_{E_{i,j}}  \cdot f_\la( X_{i,j})
$$
and put $S := (\Var T'_{\la})^{-1/2} (T'_{\la} - \E T'_{\la}) =
\sum_{i=1}^{V(\la)} (S_i - \E S_i) .$ Clearly $\Var S = \E S^2 =
1$.

We define a graph $G_\la := ({\cal V}_\la, {\cal E }_\la )$ as
follows.  The set ${\cal V}_\la$ consists of the subcubes
$Q_1,...,Q_{V(\la) }$ and  edges $(Q_i, Q_j)$ belong to ${\cal E
}_\la$ if $d(Q_i, Q_j) \leq 2 \al \la^{-1/d} \rho_\la$, where
$d(Q_i,Q_j) := \inf \{ |x - y|: x \in Q_i, y \in Q_j \}$ .
 By definition of the radius of stabilization $R(x,\la)$,
%exponential stabilization,
the value of $S_i$ is determined by the restriction of $\Po_\la$
to the $\la^{-1/d} \rho_\la$-neighborhood of the cube $Q_i$. By
the independence property of the Poisson point process, if $A_1$
and $A_2$ are disjoint collections of cubes in ${\cal V}_\la$ such
that no edge in ${\cal E }_\la$  has one endpoint in $A_1$ and one
endpoint in $A_2$, then the random variables $\{ S_{Q_i}, Q_i \in
A_1\} $ and $\{ S_{Q_j}, Q_j \in A_2 \}$ are independent.  Thus
$G_\la$ is a dependency graph for $\{S_i\}_{i=1}^{V(\la)}$.

% the random variables
%$$
% \sum_{    } \xi_\la(X_j; \P_\la  ) \cdot \1_{Q_i \cap B}(X_j) \cdot \1_{E_j }, \ \ \ \ 1 \leq i\leq V(\la),
%$$
%are independent whenever $i$ and $k$ satisfy  $d(Q_i, Q_k) \geq 2
%\al \la^{-1/d} \log \la$.

To prepare for an application of Lemma \ref{ChenShao}, we make
five observations:

(i)  $V(\la) := | {\cal V}_\la | = O(\la \rho_\la^{-d})$ as $\la \to \infty$.
% \lceil \la^{1/d}/ \rho_\la \rceil^d$.

(ii) Since the number of cubes in $Q_1,...,Q_V$ distant at most $2
%\al \la^{-1/d} \log \la$ from a given cube is bounded by $ 5^d$,
 \la^{-1/d} \rho_\la$ from a given cube is bounded by $ 5^d$,
it follows that the maximal degree $D$ satisfies $D:= D_\la \leq
5^d$.

(iii)  The definitions of $S_i$ and $\xi_{i,j}$ and Lemma
\ref{sumbound} tell us that for all $1 \leq i \leq V(\la)$ \bea
\| S_i\|_q
\leq C (\Var T'_{\la})^{-1/2}
 %\E \left[ \left( \sum_{j=1}^{\infty} |\xi_{i,j}  | \right)^3 \right]
 \left\|  \sum_{j=1}^{\infty} |\xi_{i,j} |  \right\|_q
\leq C (\Var T'_{\la})^{-1/2}
\rho_\la^{d(p+1)/p}.
\label{0511}
\eea

(iv)
We can bound $\Var[T'_\la]$ as follows. Observe that
 $T'_\la$ is the sum of $V(\la)$  random variables, which by
the case $q=2 $ of
%Jensen's inequality and
 Lemma \ref{sumbound} each have a second moment
bounded by a constant multiple of $\rho_\la^{2d(p+1)/p}$. Thus the
variance of each of the $V(\la)$ random variables is also bounded
by a constant multiple of $\rho_\la^{2d(p+1)/p}$. Moreover, the
covariance of any pair of the $V(\la)$ random variables  is zero
when the indices of the random variables correspond  to
non-adjacent cubes. For adjacent cubes, by the Cauchy-Schwarz
inequality the covariance is also bounded by a constant multiple
of $\rho_\la^{2d (p+1)/p}$.  This shows that \bea \Var[T'_\la]  =
O( \rho_\la^{d(p+2)/p} \la). \label{0511a} \eea

%(iv) $\Var [T'_{\la} ] \geq \Var [T_\la]/2$ for $\la$ large. We
(v) $\Var [T'_{\la} ]$ is close to $\Var [T_\la]$ for $\la$ large.
We require more estimates  to show this.  Note that
 $| T'_{\la} -  T_\la| = 0$  except possibly on the set
$E_\la^c $.
% which has probability less than $C \la^{-33}$.
Lemma \ref{sumbound}, along with Minkowski's inequality,
 yields the upper bound
\bea
%\E \left[ \left( \sum_{i=1}^{V(\la)}  \sum_{j = 1 }^{N_i} \vert \xi_\la
\left\|  \sum_{i=1}^{V(\la)}  \sum_{j = 1 }^{N_i} \vert \xi_\la
%((X_{i,j},U_{i,j}); \P_\la) \vert\right)^3 \right] \leq
((X_{i,j},U_{i,j}); \P_\la) \vert \1_{A_\la}(X_{i,j})
\right\|_q \leq
C V(\la) \rho_\la^{d(p+1)/p} \leq  C \la \rho_\la^{d/p}.
%C V(\la)^3 (\log \la)^{4d} \leq  C \la^4.
\label{0115}
\eea
Since $T_\la=T'_\la$ on event $E_\la$,
%Thus H\"older's inequality implies that
 the H\"older and Minkowski  inequalities
yield
\bean
 \|T_{\la} - T'_{\la}\|_2 \leq
 % \|T_{\la} -\T'_{\la}\|_q
 \|T_{\la} - T'_{\la}\|_q
 P[E_\la^c]^{(1/2)-(1/q)}
\nonumber \\
\leq
  (\|T_\la\|_q+ \|T'_\la\|_q) P[E_\la^c]^{(q-2)/(2q)}.
\eean
Hence, by  (\ref{rhodef}), (\ref{0505}), and (\ref{0115}),
\bea
 \|T_{\la} - T'_{\la}\|_2
\leq C \la \rho_\la^{d/p} (\la \tau(\rho_\la))^{(q-2)/(2q)}
 \leq C\la^{-3}
 \label{L1bounds}
\eea
which implies that
\be
 \label{L1bound}
 \E[|T'_{\la} -T_{\la} |] \leq C \la^{-3},
\ee
which we use later.
Since
$$
   \Var [T_{\la} ] = \Var[T'_\la]
+ \Var(T_\la-T'_\la) + 2 \Cov(T'_\la,T_\la - T'_\la),
$$
by (\ref{L1bounds}), (\ref{0511a}), (\ref{0513a})
 and the Cauchy-Schwarz inequality, we
obtain
\be \label{variancebounds}
  | \Var [T_{\la} ] -
\Var[T'_\la] | \leq C\la^{-2}.
\ee
\ \

Given the five observations (i)-(v), we are now ready to
 apply Lemma \ref{ChenShao} to prove Theorem \ref{mainthm1}.
By (\ref{0511a}) and (\ref{variancebounds}), $\Var[T_\la]$, as a
function of $\la$, is bounded on bounded intervals. Hence, it
suffices to prove that there exists $\la_0  \geq 2$ such that
(\ref{rates2})  holds for all $\la \geq \la_0$, since we can then
extend (\ref{rates2}) to all $\la \geq 2$ by changing $C$ if
necessary.  Trivially, (\ref{rates2}) holds for large enough $\la$
when   $\Var [T_{\la}] < 1$, and so
   without loss of generality we
  now  assume $\Var [T_{\la}] \geq 1$.
To establish the error bound (\ref{rates2}) in this case, we apply
the bound (\ref{CS1})  to $W_i := S_i - \E S_i, \ 1 \leq i \leq
V_\la,$
%with $p =3$ and
 with
$$
\theta : = C ( \Var T'_{\la} )^{-1/2}
%  ( \log \la)^{4d/3}.
 \rho_\la^{d(p+1)/p}.
$$
Our choice of $\theta$ is applicable by (\ref{0511}).
% because of observation (iii).
We clearly have $\E [W_i] = 0$ and $\E [(\sum_{i=1}^{V(\la)}
W_i)^2] = 1$. With $S = \sum_{i=1}^{V(\la)} W_i$, Lemma
\ref{ChenShao} along with observation (i) above yields
 $$
 \sup_t \left| P[S \leq t] - \Phi(t) \right|
 \leq C
  \la \rho_\la^{-d}
 ( \Var T'_{\la} )^{-q/2}  \rho_\la^{dq(p+1)/p}
$$
 \be \label{diffbound} \leq C   \la   ( \Var T_{\la} )^{-q/2}
\rho_\la^{dq},
%( \log \la)^{3d},
 \ee where the last line makes use of the fact that
 $\Var [T'_{\la}] \geq
\Var [T_\la ]/2$, which follows (for $\la$ large)  from
 (\ref{variancebounds}).

Now  if $\beta > 0$ is a constant and $Z$ any random variable then
by (\ref{diffbound}) we have for all $t \in \R$
$$
P[Z \leq t] \leq P[S \leq t + \b ]  + P[| Z - S | \geq \b] $$
$$
\leq \Phi(t + \b) + C \la
 ( \Var T_{\la} )^{-q/2} \rho_{\la}^{dq} + P[| Z - S | \geq \b]
$$
$$
\leq \Phi(t ) +  C \b + C \la   ( \Var T_{\la} )^{-q/2}
\rho_\la^{dq} + P[| Z - S | \geq \b]
$$
by the Lipschitz property of $\Phi$. Similarly for all $t \in \R$,
$$
P[Z \leq t] \geq  \Phi(t ) -  C \b - C  \la   ( \Var T_{\la}
)^{-q/2} \rho_\la^{dq} - P[| Z - S | \geq \b].
$$
In other words \be \label{bound4}  \sup_t | P[Z \leq t] - \Phi(t)
| \leq C \b + C \la   ( \Var T_{\la} )^{-q/2}  \rho_\la^{dq}  +
P[| Z - S | \geq \b]. \ee

Now by definition of $S$,
$$
|(\Var T'_{\la})^{-1/2} (T_{\la} - \E T_{\la}) - S| = | (\Var
T'_{\la})^{-1/2} \{ (T_{\la} - \E T_{\la}) - (T'_{\la} - \E
T'_{\la}) \} |
$$
$$
  \leq (\Var T'_{\la})^{-1/2} \{ | T_{\la} - T'_{\la} | + \E [ | T_{\la} - T'_{\la} | ] \}
$$
which by (\ref{L1bound}) is bounded by $C \la^{-3}$ except
possibly on the set $E_\la^c$ which has probability less than $C
\la^{-2}$ by (\ref{0505}) and (\ref{rhodef}).
 Thus by (\ref{bound4}) with  $Z = (\Var
T'_{\la})^{-1/2} (T_{\la} - \E T_{\la})$ and $\beta = C \la^{-3}$
\bea
 \sup_t | P[(\Var T'_{\la})^{-1/2} (T_{\la} -
  \E T_{\la} )  \leq t] - \Phi(t) |  \leq C   \la  ( \Var
T_{\la} )^{-q/2}
%( \log \la)^{3d}
 \rho_\la^{dq}
 +  C \la^{-2}.
%+ C \la^{-33}.
\label{0627}
\eea

Moreover, by the triangle inequality
\bea
 \sup_t \left| P[(\Var T_{\la})^{-1/2} (T_{\la} - \E T_{\la}) \leq t] - \Phi(t) \right| \leq
\nonumber \\
\leq  \sup_t \left|  P\left[ (\Var T'_{\la})^{-1/2} (T_{\la} - \E
T_{\la}) \leq t \cdot ( {\Var T_{\la} \over \Var T'_{\la} }
)^{1/2} \right] - \Phi \left(t ( {\Var T_{\la} \over \Var T'_{\la}
} )^{1/2} \right) \right|  +
\nonumber \\
 + \ \sup_t \left| \Phi \left(t ( {\Var T_{\la} \over \Var T'_{\la} } )^{1/2} \right) -
 \Phi(t)
\right|.
\label{0513}
\eea
Since for all $s \leq t$, we have $| \Phi(s) - \Phi(t) | \leq (t -
s) \max_{s \leq u \leq t} \phi(u)$ where $\phi$ denotes the
standard normal density, and since by (\ref{variancebounds}) there
is a constant $0 < C < \infty$ such that for all $\la > 0$ and all
$t \in \R$
$$
 \left|t ( {\Var T_{\la} \over \Var T'_{\la} } )^{1/2} - t \right| \leq |t|
\left|  {\Var T_{\la} \over \Var T'_{\la} }  - 1 \right|
 \leq {C|t| \over \la^2}
$$
we get
$$
 \sup_t \left| \Phi \left(t ( {\Var T_{\la} \over \Var T'_{\la} } )^{1/2} \right) -
\Phi(t) \right| \leq  C \sup_t \left( \left({ |t | \over \la^2 }
\right) \left( \max_{u \in [ t - tC/\la^2, \ t +  tC/\la^2]}
\phi(u) \right) \right) \leq {C \over \la^2}.
$$
Thus by (\ref{0627}) and (\ref{0513}),  \bea
 \sup_t| P[(\Var T_{\la})^{-1/2} (T_{\la} - \E T_{\la}) \leq t] - \Phi(t) |
\leq  C \la     ( \Var T_{\la} )^{-q/2} \rho_\la^{dq} +  C \la^{-2}.
\label{lasteq}
\eea

Finally we can deduce from (\ref{0511a}) and (\ref{variancebounds})
 that $\Var T_\la  =  O( \la \rho_\la^{d(p+2)/p} )$.
Hence, under the assumptions of Theorem \ref{mainthm1},
  in (\ref{lasteq}) the first term in the right hand
 side dominates,
thus  yielding the desired bound (\ref{rates2}), and the proof of
Theorem \ref{mainthm1} is complete.

Under the assumptions of Theorem \ref{polyth}, provided $\sigma^2
(f,\xi,\ka)>0$ the right hand side of (\ref{lasteq}) is bounded by
$C\la^{-1/2} \rho_\la^{dq}$, and since in this case we set
$\rho_\la = \la^a$ with $a$ given by (\ref{adef}), some elementary
algebra yields (\ref{rates2a}). Our assumption that $\gamma >
d(150 + 6/p)$ then yields the central limit theorem behavior
(\ref{CLTeq}), which is also trivially true in the case where
$\sigma^2 (f,\xi,\ka)=0$. This completes the proof of Theorem
\ref{polyth}. \qed

\vskip.5cm

%\Comment{added acknowledgment}

 {\bf Acknowledgments.} We began
this work while visiting the Institute for Mathematical Sciences
at the National University of Singapore, and continued it while
visiting the Isaac Newton Institute for Mathematical Sciences at
Cambridge. We thank both institutions for their hospitality.

%Mathew D. Penrose, Department of Mathematical Sciences,
%University of Bath, Bath BA2 7AY United Kingdom:
%{\texttt mathew.penrose@durham.ac.uk}

\end{document}